\documentclass[smallextended,final]{svjour3}[27.12.2011] 
\usepackage{booktabs}
\usepackage{multirow}
\usepackage{amsmath,amsfonts}
\usepackage{colortbl}
\usepackage{hyperref}
\usepackage{graphicx}
\usepackage{algorithmic, algorithm}

\definecolor{gray}{rgb}{0.80,0.80,0.80}
\definecolor{lightgray}{rgb}{0.92,0.92,0.92}

\newcommand{\bt}{{\beta}}

\DeclareMathOperator*{\rank}{rank}

\title{Multilevel quasiseparable matrices in PDE-constrained optimization}

\author{Jacek Gondzio \and Pavel Zhlobich}

\institute{School of Mathematics, The University of Edinburgh, JCMB, The King's Buildings, Edinburgh, Scotland EH9 3JZ, \email{J.Gondzio@ed.ac.uk, P.Zhlobich@ed.ac.uk}.}

\begin{document}

\titlerunning{Multilevel quasiseparable matrices in PDE-constrained optimization}

\maketitle

\begin{abstract}
Optimization problems with constraints in the form of a partial differential equation arise frequently in the process of engineering design. The discretization of PDE-constrained optimization problems results in large-scale linear systems of saddle-point type. In this paper we propose and develop a novel approach to solving such systems by exploiting so-called \emph{quasiseparable} matrices. One may think of a usual quasiseparable matrix as of a discrete analog of the Green's function of a one-dimensional differential operator. Nice feature of such matrices is that almost every algorithm which employs them has linear complexity. We extend the application of quasiseparable matrices to problems in higher dimensions. Namely, we construct a class of preconditioners which can be computed and applied at a linear computational cost. Their use with appropriate Krylov methods leads to algorithms of nearly linear complexity.


\end{abstract}

\section{Introduction}\label{sec:intro}

In this paper we introduce a new class of structured matrices that we suggest to call \emph{multilevel quasiseparable}. It generalizes well-known quasiseparable matrices to the multi-dimensional case. Under the multilevel paradigm, parameters that are used to represent a matrix of a higher hierarchy are themselves multilevel quasiseparable of a lower hierarchy. The usual one-dimensional quasiseparable matrix is the one of the lowest hierarchy.

Quasiseparable matrices found their application in orthogonal polynomials \cite{BEGOZ10}, root finding \cite{BDG04}, system theory \cite{BOZ10b} and other branches of applied mathematics and engineering. Consequently, there has been major interest in them in the scientific community, and many fast algorithms for working with them have been developed in the last decade. Due to the similarity in representations of multilevel quasiseparable and usual quasiseparable matrices, it is possible to extend the applicability of almost all the algorithms developed for quasiseparable matrices to the new class. It is a celebrated fact that operations with quasiseparable matrices can be performed in linear time with respect to their size. We show that under some conditions this is also true for multilevel quasiseparable matrices.

Natural applications of multilevel quasiseparable matrices are discretized Partial Differential Equations (PDEs) and related problems. There have already been recent attempts to use structured matrices similar to quasiseparable ones to solve discretized PDEs \cite{XCGL09,M09}. The very distinct feature of our approach is the ability not only to solve systems in linear time, but also to \emph{invert} matrices and to perform arithmetic operations with them in linear time. Therefore, we can solve, for instance, saddle-point systems with PDE matrices that former approaches did not allow.

The structure of the paper is as follows: we start with a motivation for this paper by considering a simple PDE-constrained optimization problem in Section \ref{sec:motivation}. This type of problems leads to systems in a saddle-point form. The need in efficient solvers for such systems motivated this paper. In Section \ref{sec:qs} we briefly review quasiseparable matrices and related results. Section \ref{sec:mlqs} concerns the extension of quasiseparable matrices to the multi-dimensional case. We give there a formal definition of multilevel quasiseparable matrices, introduce their parametrization and fast arithmetic with them. The last section presents results of preliminary numerical experiments that empirically confirm linear time complexity of algorithms with multilevel quasiseparable matrices.

\section{Model PDE-constrained optimization problem.}\label{sec:motivation}
As a simple model let us consider a distributed control problem which is composed of a cost functional
\begin{equation}\label{eq:cost_functional}
\min_{u,f} \frac{1}{2}\|u-\widehat{u}\|_2^2+\beta\|f\|^2
\end{equation}
to be minimized subject to a Poisson's problem with Dirichlet boundary conditions:
\begin{equation}\label{eq:constraint}
\left.\begin{array}{l}
-\nabla^2u=f\mbox{ in }\Omega,\\
u=g\mbox{ on }\partial\Omega.
\end{array}\right.
\end{equation}

The problem consists in finding a function $u$ that satisfies the PDE constraint and is as close as possible in the L2 norm sense to a given function $\widehat{u}$ (``desired state''). The right-hand side of the PDE $f$ is not fixed and can be designed in order to achieve the optimality of $u$. In general such a problem is ill-posed and thus needs Tikhonov regularization that is introduced in the form of the second term in \eqref{eq:cost_functional}.

The usual approach to solving the minimization problem \eqref{eq:cost_functional} -- \eqref{eq:constraint} numerically is to introduce a weak formulation of it and to discretize the latter by finite elements, see, for instance, \cite{RDW10,ESW05}. After these steps, the discrete analog of \eqref{eq:cost_functional} -- \eqref{eq:constraint} reads
\begin{equation}\label{eq:discrete_PDE_constrained}
\begin{aligned}
&\min_{u_h,f_h}\frac{1}{2}\|u_h-\widehat{u}\|_2^2 + \beta\|f_h\|_2^2,\quad u_h\in V_g^h, f_h\in V_0^h,\\
&\mbox{subject to } \int\limits_{\Omega}\nabla u_h\cdot\nabla v_h = \int\limits_{\Omega}v_h f,\quad \forall v_h\in V_0^h,
\end{aligned}
\end{equation}
where $V_0^h$ and $V_g^h$ are finite-dimensional vector spaces spanned by test functions. If $\{\phi_1,\ldots\phi_n\}$ and $\{\phi_1,\ldots\phi_n,\phi_{n+1},\ldots\phi_{n+\partial n}\}$ are bases in $V_0^h$ and $V_g^h$, respectively, and $u_h$ and $v_h$ are represented in them as
\[
u_h = \sum\limits_{j = 1}^{n+\partial n}u_j\phi_j,\quad f_h = \sum\limits_{j = 1}^{n}f_j\phi_j,
\]
then the matrix form of \eqref{eq:discrete_PDE_constrained} is
\begin{equation}\label{eq:PDE_constrained_matrix_form}
\begin{aligned}
&\min_{\vec{u},\vec{f}} \frac{1}{2}\vec{u}^TM\vec{u} - \vec{u}^T\vec{b} + c + \beta \vec{f}^TM\vec{f},\\
&\mbox{subject to } K\vec{u} = M\vec{f} + \vec{d},
\end{aligned}
\end{equation}
where $M_{ij}=\int_{\Omega}\phi_i\phi_j$ --- mass matrix, $K_{ij}=\int_{\Omega}\nabla\phi_i\nabla\phi_j$ --- stiffness matrix, $b_j=\int_{\Omega}\phi_j\widehat{u}$, $c=\frac{1}{2}\int_{\Omega}\widehat{u}^2$, $d_j=\sum\limits_{k=n+1}^{n+\partial n}u_k\int_{\Omega}\nabla\phi_k\nabla\phi_j$, $i,j=1,\ldots,n$.

The Lagrangian associated with problem
\eqref{eq:PDE_constrained_matrix_form} is
\[
\mathcal{L}(\vec{u},\vec{f},\vec{\lambda}) = \frac{1}{2}\vec{u}^TM\vec{u} - \vec{u}^Tb + c + \beta \vec{f}^TM\vec{f} + \vec{\lambda}^T(K\vec{u} - M\vec{f} - \vec{d}),
\]
where $\vec{\lambda}$ is a vector of Lagrange multipliers. The condition for a stationary point of $\mathcal{L}$ define $\vec{u}$, $\vec{f}$ and $\vec{\lambda}$ via the solution of the linear system
\begin{equation}\label{eq:saddle-point}
\begin{bmatrix}
2\beta M & 0 & -M \\
0 & M & K^T \\
-M & K & 0
\end{bmatrix}
\begin{bmatrix}
\vec{f} \\
\vec{u} \\
\vec{\vec{\lambda}}
\end{bmatrix}
=
\begin{bmatrix}
\vec{0} \\
\vec{b} \\
\vec{d}
\end{bmatrix}.
\end{equation}

This linear system is of saddle-point type. We refer to \cite{BGL05} for an extensive survey of numerical methods for this type of systems.

Matrix in \eqref{eq:saddle-point} is symmetric, usually very large but sparse due to finite element discretization. Thus matrix-vector multiplication is cheap and Krylov subspace methods such as Minimal Residual (MINRES) method of Paige and Saunders \cite{PS75} or Projected Preconditioned Conjugate Gradient (PPCG) method of Gould, Hribar and Nocedal \cite{GHN02} are well suited. Their convergence nevertheless depends highly on the choice of a good preconditioner \cite{RDW10}.

In this paper we propose a structured matrices approach to the solution of systems of the type described above. First, let us write the block LDU decomposition of \eqref{eq:saddle-point}:
\begin{equation}\label{eq:LDU}
\begin{gathered}
\begin{bmatrix}
2\beta M & 0 & -M \\
0 & M & K^T \\
-M & K & 0
\end{bmatrix}
=
\begin{bmatrix}
I & 0 & 0 \\
0 & I & 0 \\
\frac{-1}{2\beta}I & KM^{-1} & I
\end{bmatrix}
\cdot
\begin{bmatrix}
2\beta M & 0 & 0 \\
0 & M & 0 \\
0 & 0 & S
\end{bmatrix}
\cdot
\begin{bmatrix}
I & 0 & \frac{-1}{2\beta}I \\
0 & I & M^{-1}K^T \\
0 & 0 & I
\end{bmatrix},\\ S = -\left(\frac{1}{2\beta}M+KM^{-1}K^T\right).
\end{gathered}
\end{equation}
The hardest part in using this decomposition for solving a system is to compute the Schur complement $S$ of the (3,3) block and to solve the corresponding system with it. Since matrix $S$ is usually dense, this task is untractable if we use its entrywise representation but it is feasible if we use a proper structured representation. In Section \ref{sec:mlqs} we will show that Schur complement indeed has a structure and this structure is the one called \emph{multilevel quasiseparable}. We will give all necessary definitions and show that the use of multilevel quasiseparable matrices leads to asymptotically (with respect to the size of the system) linear storage and linear complexity algorithm for solving \eqref{eq:saddle-point}.

\section{Quasiseparable matrices.}\label{sec:qs}
Matrix of a rank much smaller than its size is a discrete analog of a separable function. More generally, we may think of a matrix with certain blocks being low rank rather than the whole matrix. This case corresponds to a function being separable in some subdomains of the whole domain. One simple example of such a function is Green's function of the Sturm--Liouville equation (this example will be considered in some more details later in this section). There is no surprise that matrices with low rank blocks found their applications in many branches of applied mathematics, especially in integral equations. There are several related classes of rank-structured matrices that one can find in the scientific literature. Among the most well-known are semiseparable \cite{GKK85}, quasiseparable \cite{EG99a}, $H$-matrices \cite{H99}, $H^2$-matrices \cite{HKS00} and mosaic-skeleton matrices \cite{T96}. In the current paper our main tool will be quasiseparable matrices \cite{VVM05} having low rank partitions in their upper and lower parts. The formal definition is given below.

\begin{definition}[Rank definition of a block quasiseparable matrix]\label{def:qsrank}
Let $A$ be a block matrix of block sizes $\{n_k\}_{k=1}^n$ then it is called block $(r^l,r^u)$-quasiseparable if
\[
\begin{gathered}
\max_{K}\,\rank A(K+1:N,1:K)\leq r^l,\quad\max_{K}\,\rank A(1:K,K+1:N)\leq r^u,\\
K=\sum\limits_{i=1}^k n_i,\quad N=\sum\limits_{i=1}^n n_i,
\end{gathered}
\]
where $r^l$ ($r^u$) is called lower (upper) order of quasiseparability.
\end{definition}

In what follows we will call \emph{block $(r^l,r^u)$-quasiseparable} matrices simply \emph{quasiseparable} for shortness.

A well-known example of a quasiseparable matrix is given by the discretized Green's function of a one-dimensional differential equation. Consider a regular inhomogeneous Sturm-Liouville boundary value problem:
\begin{equation}\label{eq:S-L}
\begin{aligned}
&(p(x)u')'-q(x)u=f(x),\\
&\begin{cases}
\alpha_1 u(a)+\beta_1 u(a) = 0,&\\
\alpha_2 u(b)+\beta_2 u(b) = 0,&\\
\end{cases}
\end{aligned}
\end{equation}
where functions $p$, $p'$, $q$ are continuous on $[a,b]$, $p(x)>0$ on $[a,b]$ and $|\alpha_i|+|\beta_i|\neq 0$ for $i=1,2$. It is a classical result in ordinary differential equations that the solution $f(x)$ of \eqref{eq:S-L} is given by
\[
u(x)=\int\limits_{a}^{b}g(x,\xi)f(\xi)\,\mathrm{d}\xi,
\]
where $g(x,\xi)$ is the so-called Green's function. In this case it has an explicit formula
\begin{equation}\label{eq:GreensF}
g(x,\xi) = \frac{1}{p(a)W(a)}=
\begin{cases}
u_1(x)u_2(\xi),&a \leq x \leq \xi,\\
u_1(\xi)u_2(x),&\xi < x \leq b,\\
\end{cases}
\end{equation}
with $W(x)$ being Wronskian and $u_1(x)$ and $u_2(x)$ being two specific linearly independent solutions of \eqref{eq:S-L}. It is easy to see from \eqref{eq:GreensF} that discretized Green's function is a quasiseparable matrix of order one. 

In order to exploit the quasiseparability of matrices in practice one must have a low-parametric representation of them. There are many alternative parametrisations of quasiseparable matrices all of which use $\mathcal{O}(N)$ parameters, where $N$ is the size of the matrix. Having such a parametrisation at hand one can write most of the algorithms, e.g., inversion, LU or QR decomposition, matrix-vector multiplication in terms of these parameters and, therefore, the complexity of these algorithms is $\mathcal{O}(N)$. In this paper we will use the so-called \emph{generator representation} (Definition \ref{def:qsgen} below) proposed by Eidelman and Gohberg \cite{EG99a}. There are several alternative parametrisations such as Givens-weight \cite{DV07} and Hierarchical Semiseparable (HSS) \cite{CGP06}.

\begin{definition}[Generator definition of a block quasiseparable matrix]\label{def:qsgen}
Let $A$ be a block matrix of block sizes $\{n_k\}_{k=1}^n$ then it is called block $(r^l,r^u)$-quasiseparable if it can be represented in the form
\begin{equation}\label{eq:qs_matrix}
\begin{bmatrix}
d_1 & g_1h_2 & g_1b_2h_3 & \cdots & g_1b_2\dots b_{n-1}h_n \\
p_2q_1 & d_2 & g_2h_3 & \cdots & g_2b_3\dots b_{n-1}h_n \\
p_3a_2q_1 & p_3q_2 & d_3 & \cdots & g_3b_4\dots b_{n-1}h_n \\
\vdots & \vdots & \vdots & \ddots & \vdots \\
p_na_{n-1}\dots a_2q_1 & p_na_{n-1}\dots a_3q_2 & p_na_{n-1}\dots
a_4q_3 & \cdots & d_n
\end{bmatrix},
\end{equation}
where parameters (called generators) $\{d_k$, $q_k$, $a_k$, $p_k$, $g_k$, $b_k$, $h_k\}$ are matrices of sizes as in the table below.

\begin{table}[ht]
\centering
\caption{Sizes of generators of a block quasiseparable matrix.}\label{tbl:uhgens}
\begin{tabular}[c]{l>{$}c<{$}>{$}c<{$}>{$}c<{$}>{$}c<{$}>{$}c<{$}>{$}c<{$}>{$}c<{$}}
& d_k & q_k & a_k & p_k & g_k & b_k & h_k\\
\midrule
\# of rows & n_k & r^l_k & r^l_k & n_{k} & n_k & r^u_{k-1} & r^u_{k-1} \\
\midrule
\# of cols & n_k & n_k & r^l_{k-1} & r^l_{k-1} & r^u_k & r^u_k & n^u_k \\
\bottomrule
\end{tabular}
\end{table}

Orders of quasiseparability $r^l$ and $r^u$ are maxima over the corresponding sizes of generators:
\[
r^l = \max_{1\leq k\leq n-1} r^l_k,\quad r^u = \max_{1\leq k\leq n-1} r^u_k.
\] 
\end{definition}

\begin{remark}
Generators are not unique, there are infinitely many ways to represent the same quasiseparable matrix with different generators. For the relation between different sets of generators see \cite{EG05}.
\end{remark}

\begin{remark}
It is obvious that any scalar quasiseparable matrix can be converted to the block form by simple aggregation of its generators.
\end{remark}

\begin{remark}
Sizes $r^l_k$ and $r^u_k$ of generators are directly related to ranks of submatrices in the lower and upper parts, respectively. Namely, if we let $K$ to be the block index: $K=\sum\limits_{i=1}^k n_i$, then
\begin{equation}\label{eq:subm_ranks}
\rank A(K+1:N,1:K) \leq r^l_k,\quad\rank A(1:K,K+1:N) \leq r^u_k.
\end{equation}
Moreover, for any quasiseparable matrix there exists a set of generators with sizes $r^l_k$ and $r^u_k$ that satisfy \eqref{eq:subm_ranks} with exact equalities (such generators are called \emph{minimal}). For their existence and construction see \cite{EG05}.
\end{remark}

One remarkable property of quasiseparable matrices is that this class is closed under inversion \cite{EG99a,SN04}. For instance, inverse of a banded matrix is not banded but both matrices are quasiseparable. Due to the low-parametric representation in Definition \ref{def:qsgen} most of the algorithms with quasiseparable matrices have linear complexities. Table \ref{tbl:alg} lists these algorithms along with the corresponding references.

\begin{table}[ht]
\centering
\caption{$\mathcal{O}(n)$ algorithms for quasiseparable matrices}\label{tbl:alg}
\begin{tabular}{ccccccc}
$A * v$ & LU & QR & $A^{-1}$ & $A*B$ & $A\pm B$ & $\|A\|_F$ \\
\midrule
\cite{EG99b} & Theorem \ref{thm:LU} & \cite{EG02} & \cite{EG99a} & \cite{EG99a} & -- & \cite{EG05} \\
\bottomrule
\end{tabular}
\end{table}

The ability to solve systems with quasiseparable matrices in $\mathcal{O}(n)$ operations is essential for the purpose of this paper. One of the ways to do it is through the use of fast LU decomposition. We next derive LU decomposition algorithm for a general block quasiseparable matrix in terms of the generators it is represented with. A restricted version of this algorithm applicable to diagonal plus semiseparable matrices (a subclass of quasiseparable matrices) was derived in \cite{GKK85} and some formulae of the current algorithm are similar to the ones in the inversion algorithm given in \cite{EG99a}. Still, to the best of our knowledge, current LU decomposition algorithm has never been published and may be useful to those who need a fast system solver for quasiseparable matrices. The complexity of the algorithm is $\mathcal{O}(N)$ and it is valid in the \emph{strongly regular} case (i.e. \emph{all of the principal leading minors are non-vanishing}). First, let us note that quasiseparable structure of the original matrix implies the quasiseparable structure of $L$ and $U$ factors. The theorem below makes this statement precise and, in addition, relates generators of an original matrix to generators of the factors.

\begin{theorem}\label{thm:LU}
Let $A$ be a strongly regular $N\times N$ block $(r^l,r^u)$-quasiseparable matrix given by generators $\{d_k,q_k,a_k,p_k,g_k,b_k,h_k\}$ as in \eqref{eq:qs_matrix}. Let $A=LU$ be its block LU decomposition of the same block sizes. Then
\begin{enumerate}
\item[(i)] Factors $L$ and $U$ are $(r^l,0)$-- and $(0,r^u)$-quasiseparable. Moreover, $r^l_k(L)=r^l_k(A)$ and $r^u_k(U)=r^u_k(A)$ for all $k=1,\ldots,n-1$.
\item[(ii)] $L$ and $U$ are parametrized by the generators
$\{I_k,\widetilde{q}_k,a_k,p_k,0,0,0\}$ and $\{\widetilde{d}_k$, $0$, $0$, $0$, $\widetilde{g}_k$, $b_k$, $h_k\}$, where $I_k$ are identity matrices of sizes $n_k\times n_k$ and new parameters $\widetilde{q}_k$, $\widetilde{d}_k$ and $\widetilde{g}_k$ can be computed using the following algorithm:
\begin{algorithm}[ht]
\caption{Fast quasiseparable LU decomposition.}\label{alg:LU}
\begin{algorithmic}[1]
\REQUIRE $d_k,q_k,a_k,p_k,g_k,b_k,h_k$
\STATE $\widetilde{d}_1=d_1,\quad \widetilde{q}_1=q_1\widetilde{d}_1^{-1},\quad \widetilde{g}_1=g_1,\quad f_1 = \widetilde{q}_1\widetilde{g}_1$
\FOR {$k=2$ \TO $n-1$}
\STATE $\widetilde{d}_k=d_k - p_kf_{k-1}h_k$
\STATE $\widetilde{q}_k=(q_k - a_kf_{k-1}h_k)\widetilde{d}_k^{-1}$
\STATE $\widetilde{g}_k=g_k - p_kf_{k-1}b_k$
\STATE $f_k = a_kf_{k-1}b_k+\widetilde{q}_k\widetilde{g}_k$
\ENDFOR
\STATE $\widetilde{d}_n=d_n - p_nf_{n-1}h_n.$
\ENSURE $\widetilde{d}_k,\widetilde{q}_k,\widetilde{g}_k$
\end{algorithmic}
\end{algorithm}
\end{enumerate}
\end{theorem}
\begin{proof}
Statement (i) of the theorem follows from statement (ii), so we only need to prove the latter part. 

Denote, as usual, by $K$ the block index: $K=\sum\limits_{i=1}^k n_i$ and note that quasiseparable representation \eqref{eq:qs_matrix} implies the following recurrence relation between the blocks of $A$:
\begin{equation}\label{eq:block_recurrences}
\begin{aligned}
&A=
\left[\begin{array}{c|c}
A(1:K,1:K) & G_kH_{k+1} \\
\hline
P_{k+1}Q_k & A(K+1:N,K+1:N)
\end{array}\right];\\
&Q_1 = q_1,\quad Q_k = [a_kQ_{k-1}\;\;q_k],\quad k = 2,\ldots n-1;\\
&P_n = p_n,\quad P_k = [p_k\;;\;P_{k+1}a_k],\quad k = n-1,\ldots 2;\\
&G_1 = g_1,\quad G_k = [G_{k-1}b_k\;;\;g_k],\quad k = 2,\ldots n-1;\\
&H_n = h_n,\quad H_k = [h_k\;\;b_kH_{k+1}],\quad k = n-1,\ldots 2.
\end{aligned}
\end{equation}

The proof will be constructed by induction. We will show that for each $k$
\begin{equation}\label{eq:blockLU}
\left[\begin{array}{c|c}
A_{11}^k & G_kH_{k+1} \\
\hline
P_{k+1}Q_k & \star
\end{array}\right]=
\left[\begin{array}{c|c}
L_{11}^k & 0 \\
\hline
P_{k+1}\widetilde{Q}_k & \star
\end{array}\right]
\cdot
\left[\begin{array}{c|c}
U_{11}^k & \widetilde{G}_kH_{k+1} \\
\hline
0 & \star
\end{array}\right].
\end{equation}

For $k=1$ we get from \eqref{eq:blockLU}:
\[
\begin{aligned}
&d_1 = \widetilde{d}_1,\\
&P_2Q_1=P_2\widetilde{Q}_1\widetilde{d}_1,\\
&G_1H_2 = \widetilde{G}_1H_2,
\end{aligned}\quad\Longleftarrow\quad
\begin{aligned}
&\widetilde{d}_1 = d_1,\\
&\widetilde{q}_1 = q_1\widetilde{d}_1^{-1},\\
&\widetilde{g}_1= g_1.
&\end{aligned}
\]

Let us introduce an auxiliary parameter $f_k=\widetilde{Q}_k\widetilde{G}_k$. It is easy to show by using \eqref{eq:block_recurrences} that $f_k$ satisfies the recurrence relation 
\begin{equation}\label{eq:f_k}
f_1 = \widetilde{q}_1\widetilde{g}_1,\quad f_k = a_kf_{k-1}b_k+\widetilde{q}_k\widetilde{g}_k.
\end{equation}

Assume that \eqref{eq:blockLU} holds for some fixed $k$, then it holds for $k+1$ if
\begin{eqnarray}
d_{k+1} = [p_{k+1}\widetilde{Q}_k \;\; 1]\cdot[\widetilde{G}_kh_{k+1} \;;\; \widetilde{d}_{k+1}],\label{eq1}\\
P_{k+2}q_{k+1} = P_{k+2}\widetilde{Q}_{k+1}\cdot[\widetilde{G}_kh_{k+1} \;;\; \widetilde{d}_{k+1}],\label{eq2}\\
g_{k+1}H_{k+2} = [p_{k+1}\widetilde{Q}_k \;\; 1]\cdot\widetilde{G}_{k+1}H_{k+2}.\label{eq3}
\end{eqnarray}

The first equality simplifies
\[
d_{k+1} = p_{k+1}\widetilde{Q}_k\widetilde{G}_kh_{k+1}+\widetilde{d}_{k+1} = p_{k+1}f_kh_{k+1}+\widetilde{d}_{k+1}.
\]
The second equality \eqref{eq2} holds if
\begin{multline*}
q_{k+1} = \widetilde{Q}_{k+1}\cdot[\widetilde{G}_kh_{k+1} \;;\; \widetilde{d}_{k+1}] = \\ = [a_{k+1}\widetilde{Q}_{k}\;\;\widetilde{q}_{k+1}]\cdot[\widetilde{G}_kh_{k+1} \;;\; \widetilde{d}_{k+1}] =
a_{k+1}f_kh_{k+1}+\widetilde{q}_{k+1}\widetilde{d}_{k+1}.
\end{multline*}
Finally, the last equality \eqref{eq3} is true if
\[
g_{k+1}=[p_{k+1}\widetilde{Q}_k \;\; 1]\cdot\widetilde{G}_{k+1} = 
[p_{k+1}\widetilde{Q}_k \;\; 1]\cdot[\widetilde{G}_{k}b_{k+1}\;;\;\widetilde{g}_{k+1}] = 
p_{k+1}f_kb_{k+1}+\widetilde{g}_{k+1}.
\]

Matrix $\widetilde{d}_{k+1}$ is invertible because, by our assumption, matrix $A$ is strongly regular. Hence, we conclude that \eqref{eq:blockLU} is true also for index $k+1$ if generators $\widetilde{q}_{k+1}$, $\widetilde{d}_{k+1}$ and $\widetilde{g}_{k+1}$ are those computed in Algorithm \ref{alg:LU}. The assertion of the theorem holds by induction.\qed
\end{proof}

\section{Multilevel quasiseparable matrices.}\label{sec:mlqs}
In the previous section we have emphasized the relation between second-order ordinary differential equations and quasiseparable matrices. Let us now look at the problem in two dimensions that extends \eqref{eq:S-L}:
\[
\frac{\partial}{\partial x}\left(p(x,y)\frac{\partial}{\partial x}u(x,y)\right)+
\frac{\partial}{\partial y}\left(q(x,y)\frac{\partial}{\partial y}u(x,y)\right)-
r(x,y)u(x,y)=f(x,y)
\]
for $(x,y)\in\Omega$, where $\Omega=[0,1]\times[0,1]$ with homogeneous Dirichlet boundary conditions. The standard five-point or nine-point discretization of this problem on a $n\times m$ uniform grid leads to a system of linear algebraic equations of the form
\begin{equation}\label{eq:BlockMatrix}
\begin{bmatrix}
A_1 & B_1 \\
C_1 & A_2 & B_2 \\
& C_2 & \ddots & \ddots \\
& & \ddots & \ddots & B_{m-1} \\
& & & C_{m-1} & A_m
\end{bmatrix}
\begin{bmatrix}
u_1 \\
u_2 \\
\vdots \\
\vdots \\
u_m
\end{bmatrix}
=
\begin{bmatrix}
f_1 \\
f_2 \\
\vdots \\
\vdots \\
f_m
\end{bmatrix},
\end{equation}
where we have assumed that $u_i$ are the discretized unknowns along the $i$'th column of the $n\times m$ grid. In this case each of $A_i$, $B_i$, and $C_i$ is an $n\times n$ matrix and the whole block matrix has size $nm\times nm$. Furthermore, each of the $A_i$, $B_i$, and $C_i$ is a tridiagonal matrix.

One way to solve system \eqref{eq:BlockMatrix} is to do block LU factorization assuming that it exists. When we eliminate the block entry $C_1$ we get in the position occupied by $A_2$ the new block 
\[
S_2 = A_2 - C_1A_1^{-1}B_1.
\]
Observe that even though all the individual matrices on the right-hand side of the expression are tridiagonal, $A_1^{-1}$ is not, and hence $S_2$ is a dense (non-sparse) matrix. At the next step of Gaussian elimination we would use $S_1$ as the pivot block to eliminate $C_2$. Now in the position occupied by the block $A_3$ we would get the matrix $S_3=A_3-C_2S_2^{-1}B_2$. Again, since $S_2$ is a dense matrix, in general $S_2^{-1}$ will be dense, and therefore $S_3$ will also be a dense matrix. What this implies is that during LU factorization of the sparse matrix, we will produce fill-in quickly that causes us to compute inverses (and hence LU factorizations) of dense $n\times n$ matrices. If we assume that these dense matrices have no structure, then we would need $\mathcal{O}(n^3)$ flops for that operation alone. Therefore, it follows that one would require at least $\mathcal{O}(mn^3)$ flops to compute block LU factorization of the system matrix. 

\begin{figure}
\centering
\includegraphics{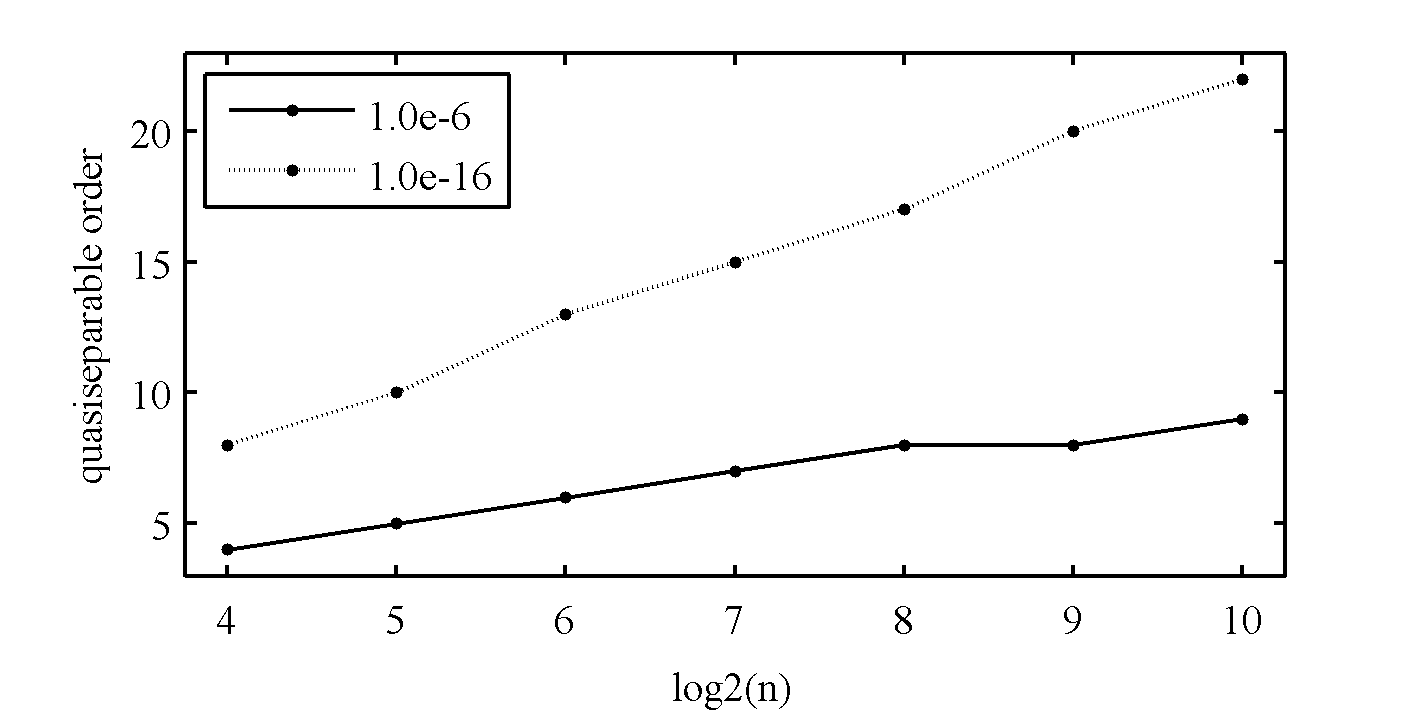}
\caption{Growth of $\varepsilon$-rank (quasiseparable order) with size for different $\varepsilon$. }\label{fig:order_vs_size}
\end{figure}

At first glance it seems that block LU factorization is not a wise approach as it does not use any kind of fill-in minimisation reordering. However, there is a rich structure in successive Schur complements $S_k$ that can be exploited to speed-up computations. In fact, it has been conjectured that if one looks at the off-diagonal blocks of these matrices ($S_2$, $S_3$, etc.), then their $\varepsilon$-rank (number of singular values not greater than $\varepsilon$) is going to be small. This conjecture
has been justified by the fact that, for example, $S_k$ can be viewed approximately (especially in the limit as $n$ becomes large) as a sub-block of the discretized Green's function of the original PDE. It is known from the theory of elliptic PDEs (see, for instance, \cite{F95}) that under some mild constraints the Green's function is smooth away from the diagonal singularity. This, in turn, implies that the numerical ranks of the off-diagonal blocks of $S_1$ are small. There have been several recent attempts to quantify this observation. It has been proved in \cite{CDGS10} that $\varepsilon$-rank of the off-diagonal blocks of Schur complements in the LU decomposition of the regular 2D-Laplacian are bounded by
\[
1+8\ln^4\left(\frac{18}{\varepsilon}\right).
\] 
This bound is not effective for any reasonable size of the problem because of the 4'th power in the logarithm (for instance for $\varepsilon=$1.0e-6 it gives 6.2e+5). However, numerical experiments with Laplacian discretized by the usual five-point stencil suggest much lower bound, see Figure \ref{fig:order_vs_size}.

Careful exploitation of the low-rank structure of Schur complements reduces the complexity of block LU factorization algorithm to linear $\mathcal{O}(mn)$. Indeed, as it has been mentioned in the previous section, all the operations with quasiseparable matrices that are needed at each step of the algorithm have $\mathcal{O}(n)$ complexity and there are $m$ steps altogether. First algorithms of this type were recently developed in \cite{XCGL09,M09}. These algorithms, although efficient in the case of simple PDEs, are not applicable to PDE-constrained optimization problems, such as the one described in Section \ref{sec:motivation}. Saddle-point matrices arising there are dense block matrices with blocks themselves being discretized PDE matrices. To use block LU factorization in this case one would need to invert PDE matrices of type \eqref{eq:BlockMatrix}. We next show that it is possible by defining the right structured representation of inverses of PDE matrices.

\begin{definition}[Multilevel quasiseparable matrix]\label{def:mlqs}
Matrix $A$ is called $d$-level quasiseparable if it admits the representation
\[
\begin{bmatrix}
D_1 & G_1H_2 & G_1B_2H_3 & \cdots & G_1B_2\dots B_{n-1}H_n \\
P_2Q_1 & D_2 & G_2H_3 & \cdots & G_2B_3\dots B_{n-1}H_n \\
P_3A_2Q_1 & P_3Q_2 & D_3 & \cdots & G_3B_4\dots B_{n-1}H_n \\
\vdots & \vdots & \vdots & \ddots & \vdots \\
P_nA_{n-1}\dots A_2Q_1 & P_nA_{n-1}\dots A_3Q_2 & P_nA_{n-1}\dots
A_4Q_3 & \cdots & D_n
\end{bmatrix},
\]
where each of the parameters $\{Q_k$, $A_k$, $P_k$, $D_k$, $G_k$, $B_k$, $H_k\}$ is a block matrix of blocks being $(d-1)$-level quasiseparable matrices. $1$-level quasiseparable matrix is a usual quasiseparable matrix, see Definition \ref{def:qsrank}.
\end{definition}

Let us give some simple examples of multilevel quasiseparable matrices to justify the usefulness of the proposed structure.

\begin{example}[Block banded matrices with banded blocks]\label{ex:Laplace}
Any banded matrix is automatically quasiseparable. Similarly, block banded matrix is 2-level quasiseparable if blocks are themselves banded matrices. As an example, consider 2D Laplacian on the square discretized by Q1 finite elements:
\begin{equation}\label{eq:Laplace}
K=
\begin{bmatrix}
A & B \\
B & \ddots & \ddots \\
& \ddots & \ddots & B \\
& & B & A
\end{bmatrix},\;
A=
\frac{1}{3}
\begin{bmatrix}
-8 & 1 \\
1 & \ddots & \ddots \\
& \ddots & \ddots & 1 \\
& & 1 & -8
\end{bmatrix},\;
B=
\frac{1}{3}
\begin{bmatrix}
1 & 1 \\
1 & \ddots & \ddots \\
& \ddots & \ddots & 1 \\
& & 1 & 1
\end{bmatrix}.
\end{equation}
\end{example}

\begin{example}[Tensor product of quasiseparable matrices]\label{ex:tensor_product} Matrix $A\otimes B$ is 2-level quasiseparable if $A$ and $B$ are 1-level quasiseparable. Hence, 2D mass matrix $M$ on the square discretized by Q1 finite elements:
\begin{equation}\label{eq:mass}
M=T\otimes T,\quad
T=
\frac{1}{6}
\begin{bmatrix}
4 & 1 \\
1 & \ddots & \ddots \\
& \ddots & \ddots & 1 \\
& & 1 & 4
\end{bmatrix}
\end{equation}
is 2-level quasiseparable. Its inverse $M^{-1}$ is also 2-level quasiseparable due to the tensor identity $(A\otimes B)^{-1}=A^{-1}\otimes B^{-1}$ and properties of 1-level quasiseparable matrices.
\end{example}

Quasiseparable orders of multilevel quasiseparable matrices may differ at each level. For instance, it follows immediately from the properties of quasiseparable matrices that inverse of a 2-level quasiseparable matrix has a representation as in Definition \ref{def:mlqs} for some generators $\{Q_k$, $A_k$, $P_k$, $D_k$, $G_k$, $B_k$, $H_k\}$. However, there is no guarantee that entries of this generators are quasiseparable matrices of the same order as for the original matrix. 2D Laplace matrix in Example \ref{ex:Laplace} is 2-level quasiseparable with generators of quasiseparable order one but as Figure \ref{fig:order_vs_size} illustrates, generators of LU factors (and, therefore, inverse) have orders larger than one. However, if it is true that numerical orders (ranks) do not differ much from the original orders we may say that inverses of multilevel quasiseparable matrices retain the structure approximately. This observation, in its turn, leads to algorithms of linear complexity.

To compute LU factorization of a 2-level quasiseparable matrix we may use Algorithm \ref{alg:LU}. It consists of approximately $12n$ arithmetic operations with 1-level quasiseparable matrices. Sum or difference of rank-$r_1$ and rank-$r_2$ matrices may have rank $r_1+r_2$. Similarly, any arithmetic operation (except inversion) performed on quasiseparable matrices may add their orders in the worst case. Therefore, to use multilevel quasiseparable structure efficiently we need some kind of order compression algorithm. Such algorithm was proposed in \cite{EG05}, we list it below for the convenience of the reader. This algorithm constitutes the essential part of the proposed method.
 
\begin{algorithm}[ht]
\caption{Quasiseparable order reduction.}\label{alg:reduction}
\begin{algorithmic}[1]
\REQUIRE $q_k$, $a_k$, $p_k$ of sizes $r_k\times n_k$, $r_k\times r_{k-1}$ and $n_k\times r_{k-1}$, respectively.
\STATE Using LQ or SVD decomposition determine matrices $L_1$ and $q_1'$ of sizes $r_1\times r_1'$ and $r_1'\times n_1$, respectively, such that $q_1=L_1q_1'$, where $q_1'(q_1')^*=I_{r_1'}$ and $r_1'=\rank q_1$.
\FOR {$k=2$ \TO $n-1$}
\STATE Using LQ or SVD decomposition determine matrices $L_k$ and $V_k$ of sizes $r_k\times r_k'$ and $r_k'\times (r_{k-1}'+n_k)$, respectively, such that $[a_kL_{k-1}\;q_k]=L_kV_k$, where $V_k'(V_k')^*=I_{r_k'}$ and $r_k'=\rank [a_kL_{k-1}\;q_k]$.
\STATE Split $V_k$ into the new generators $a_k'$, $q_k'$ of sizes $r_k'\times r_{k-1}'$, $r_k'\times n_k$, respectively: $V_k=[a_k'\; q_k']$.
\STATE $p_k'=p_kL_{k-1}$
\ENDFOR
\STATE $p_n'=p_nL_{n-1}$
\STATE Using QR or SVD decomposition determine matrices $p_n''$ and $S_{n-1}$ of sizes $n_n\times r_{n-1}''$ and $r_{n-1}''\times r_{n-1}'$, respectively, such that $p_n'=p_n''S_{n-1}$, where $(p_n'')^*p_n''=I_{r_{n-1}''}$ and $r_{n-1}''=\rank p_n'$.
\FOR {$k=n-1$ \TO $2$}
\STATE Using QR or SVD decomposition determine matrices $U_k$ and $S_{k-1}$ of sizes $(n_k+r_k'')\times r_{k-1}''$ and $r_{k-1}''\times r_{k-1}'$, respectively, such that $[p_k'\,;\,S_ka_k']=U_kS_{k-1}$, where $(U_k)^*U_k=I_{r_{k-1}''}$ and $r_{k-1}''=\rank p_k'$.
\STATE Split $U_k$ into the new generators $p_k''$, $a_k''$ of sizes $n_k\times r_{k-1}''$, $r_k''\times r_{k-1}''$, respectively: $U_k=[p_k';\,;\,a_k'']$.
\STATE $q_k''=S_kq_k'$
\ENDFOR
\STATE $q_1''=S_1q_1'$
\ENSURE New generators $q_k''$, $a_k''$, $p_k''$ of minimal possible sizes $r_k''\times n_k$, $r_k''\times r_{k-1}''$ and $n_k\times r_{k-1}''$, respectively.
\end{algorithmic}
\end{algorithm}

In exact arithmetic, new sizes $r_k''$ of generators match ranks of the corresponding submatrices in the quasiseparable matrix (see \cite{EG05} for the proof). In floating point arithmetic we may use truncated SVD or rank-revealing LQ/QR to decrease sizes of generators to the value of the $\varepsilon$-rank of the corresponding submatrix. When $\varepsilon$-ranks are small, complexity of Algorithm \ref{alg:reduction} is $\mathcal{O}(n)$ and, hence, the following fact takes place.

\emph{Complexity of algorithms listed in Table \ref{tbl:alg} with 2-level quasiseparable matrices (Definition \ref{def:mlqs}) is linear in the size if quasiseparable orders of 1-level quasiseparable matrices stay small during computations}.

In the next section we will demonstrate practical properties of multilevel quasiseparable matrices approach applied to solving PDEs and PDE-constrained optimization problems.

\section{Numerical results.}\label{sec:results}

We illustrate our new method with two different examples: 2D Laplace's equation with Dirichlet boundary conditions and distributed control problem with a constraint in the form of 2D Poisson's equation. We show that the proposed method can be used either directly or as a preconditioner for an iterative method depending on the level of order reduction used in Algorithm \ref{alg:reduction}. Although we consider symmetric problems only, our approach does not require symmetry and only depends on low rank properties of operators. Therefore, it can also be applied to convection--diffusion and Helmholtz problems. We do not consider 3D problems as it is not clear yet how to adapt SVD and LQ/QR decomposition used in Algorithm \ref{alg:reduction} to the block case.

Our software is written in C++ and is compiled with GCC compiler v.4.2.1. It is freely available for download at \url{http://www.maths.ed.ac.uk/ERGO/software.html}. All tests were done on a 2.4 GHz Intel Core 2 Duo with 4 Gb of RAM.

\begin{example}\label{ex:Laplace_equation}
Let $\Omega=[0,1]^2$ and consider the problem
\begin{equation}\label{eq:Laplace_equation}
\begin{array}{ll}
-\nabla^2u = 0,&\mbox{in }\Omega,\\
u=\widehat{u}|_{\partial\Omega},&\mbox{on }\partial\Omega,
\end{array}
\end{equation}
where
\[
\widehat{u}=
\begin{cases}
sin(2\pi y) &\mbox{if } x = 0,\\
-sin(2\pi y) &\mbox{if } x = 1,\\
0 & \mbox{otherwise}.
\end{cases}
\]
\end{example}

Let us introduce a square reqular grid on $[0,1]^2$ of mesh size $h=1/(n+1)$ and discretize equation \eqref{eq:Laplace_equation} on this grid by $Q1$ finite elements. Then equation \eqref{eq:Laplace_equation} reduces to a matrix equation with matrix $K$ in \eqref{eq:Laplace}. This matrix, as it was noted before, is 2-level quasiseparable and, therefore, we can apply fast LDU decomposition Algorithm \ref{alg:LU} with fast order-reduction Algorithm \ref{alg:reduction} at each step. Maximal quasiseparable order can be chosen adaptively in SVD but we set it manually to 4 and 8 for simplicity. Larger order corresponds to better approximation and, hence, more accurate result although it makes algorithm slower. Results of the computational tests for different mesh sizes are presented in Tables \ref{tbl:Laplace_direct_4} and \ref{tbl:Laplace_direct_8}. Note, that our solver was implemented for unsymmetric problems and did not use symmetry of the matrix in the current example. Exploiting the symmetry would roughly halve the time and memory used.

\begin{table}[ht]
\centering
\caption{Time and memory used by direct quasiseparable LDU decomposition based solver applied to the problem in Example \ref{ex:Laplace_equation}, quasiseparable order (maximal off-diagonal rank) is set to \textbf{4} during computations.}\label{tbl:Laplace_direct_4}
\begin{tabular}[c]{rrrrr}
$n$ & LDU & mem, Mb & solve & $\frac{\|Au-f\|_2}{\|f\|_2}$ \\
\midrule
$2^{12}$ & 0.11 & 3 & 0.003 & 8.22e-05 \\
$2^{14}$ & 0.51 & 12 & 0.011 & 1.85e-04 \\
$2^{16}$ & 2.18 & 51 & 0.043 & 3.93e-04 \\
$2^{18}$ & 9.00 & 210 & 0.169 & 6.91e-04 \\
$2^{20}$ & 36.88 & 847 & 0.677 & 8.81e-04 \\
\bottomrule
\end{tabular}
\end{table}

\begin{table}[ht]
\centering
\caption{Time and memory used by direct quasiseparable LDU decomposition based solver applied to the problem in Example \ref{ex:Laplace_equation}, quasiseparable order (maximal off-diagonal rank) is set to \textbf{8} during computations.}\label{tbl:Laplace_direct_8}
\begin{tabular}[c]{rrrrr}
$n$ & LDU & mem, Mb & solve & $\frac{\|Au-f\|_2}{\|f\|_2}$ \\
\midrule
$2^{12}$ & 0.19 & 4 & 0.003 & 3.31e-09 \\
$2^{14}$ & 0.91 & 19 & 0.013 & 6.19e-08 \\
$2^{16}$ & 4.09 & 83 & 0.052 & 5.72e-07 \\
$2^{18}$ & 17.44 & 342 & 0.209 & 2.33e-06 \\
$2^{20}$ & 72.83 & 1388 & 0.852 & 5.41e-06 \\
\bottomrule
\end{tabular}
\end{table}

The analysis of results collected in Tables \ref{tbl:Laplace_direct_4} and \ref{tbl:Laplace_direct_8} reveals the linear dependence of time and memory used by the algorithm on the size of the problem. The computational experience supports our claim about \emph{asymptotically linear complexity of algorithms with multilevel quasiseparable matrices}. Note also the linear dependence of time and the sub-linear dependence of memory on the maximal quasiseparable order used during computations.

If we truncate the order of quasiseparable matrices in the quasiseparable LDU decomposition to a very small number, the decomposition becomes less accurate. However, the algorithm becomes much faster and it is tempting to try this approximate LDU decomposition as a preconditioner in the preconditioned conjugate gradient method. Table \ref{tbl:Laplace_pcg} bellow illustrates the numerical performance of this approach. Results indicate that the PCG method preconditioned with  the inexact LDU decomposition is twice as fast as more accurate LDU decomposition alone. It is enough to truncate quasiseparable order to 2--4 to force PCG to converge in less than 10 iterations.

\begin{table}[ht]
\centering
\caption{Time and number of iterations (in brackets) used by PCG when applied to the problem in Example \ref{ex:Laplace_equation}, preconditioner is based on approximate quasiseparable LDU decomposition with order tolerance $r$, system is solved to the 1.0e-08 accuracy.}\label{tbl:Laplace_pcg}
\begin{tabular}[c]{rrrrr}
$n$ & $r$ & LDU & PCG & total time \\
\midrule
\multirow{2}{*}{$2^{12}$} 
 & 1 & 0.05 & 0.0359 (9) & 0.09 \\
 & 2 & 0.07 & 0.0252 (6) & 0.09 \\
\midrule
\multirow{2}{*}{$2^{14}$} 
 & 1 & 0.21 & 0.212 (14) & 0.42 \\
 & 2 & 0.30 & 0.144 (9) & 0.45 \\
\midrule
\multirow{2}{*}{$2^{16}$} 
 & 3 & 1.76 & 0.482 (7) & 2.24 \\
 & 4 & 2.18 & 0.282 (4) & 2.47 \\
\midrule
\multirow{2}{*}{$2^{18}$} 
 & 3 & 7.17 & 2.96 (11) & 10.13 \\
 & 4 & 9.05 & 1.98 (7) & 11.03 \\
\midrule
\multirow{2}{*}{$2^{20}$}
 & 4 & 37.07 & 10.1 (9) & 47.22 \\
 & 5 & 45.63 & 8.19 (7) & 53.83 \\
\bottomrule
\end{tabular}
\end{table}

Solution of problem in Example \ref{ex:Laplace_equation} obtained with PCG method with quasiseparable preconditioner is given in Figure \ref{fig:uDirect} below.

\begin{figure}
\centering
\includegraphics[scale = 0.6]{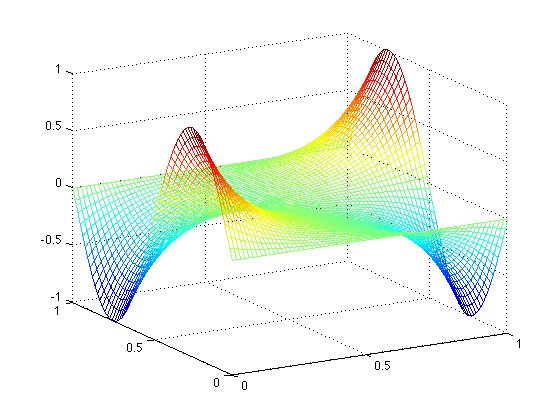}
\caption{Solution of problem in Example \ref{ex:Laplace_equation}.}
\label{fig:uDirect}
\end{figure}

\begin{example}\label{ex:Poisson_control}
Let $\Omega=[0,1]^2$ and consider the problem
\begin{equation}\label{eq:Poisson_control}
\begin{aligned}
&\min_{u,f}\frac{1}{2}\|u-\widehat{u}\|_{L^2(\Omega)}^2+\bt\|f\|_{L^2(\Omega)}^2,\\
&
\begin{array}{lll}
s.t.\;\;\;\; & -\nabla^2u = f, & \mbox{in }\Omega,\\
& u=\widehat{u}|_{\partial\Omega}, & \mbox{on }\partial\Omega,
\end{array}
\end{aligned}
\end{equation}
where
\[
\widehat{u}=
\begin{cases}
(2x-1)^2(2y-1)^2 &\mbox{if } (x,y)\in[0,\frac{1}{2}]^2,\\
0 & \mbox{otherwise}.
\end{cases}
\]
\end{example}

As it was already discussed in Section \ref{sec:motivation}, problem \eqref{eq:Poisson_control} can be solved by the discretization by Q1 finite elements. Discretized system of equation, in this case, becomes
\begin{equation}\label{eq:saddle-point2}
\begin{bmatrix}
2\beta M & 0 & -M \\
0 & M & K^T \\
-M & K & 0
\end{bmatrix}
\begin{bmatrix}
\vec{f} \\
\vec{u} \\
\vec{\lambda}
\end{bmatrix}
=
\begin{bmatrix}
\vec{0} \\
\vec{b} \\
\vec{d}
\end{bmatrix}.
\end{equation}
After we eliminate variables $\vec{f}$ and $\vec{u}$ in equation \eqref{eq:saddle-point2} it reduces to 
\begin{equation}\label{eq:normal-equation}
\left(\frac{1}{2\beta}M+KM^{-1}K^T\right)\vec{\lambda}=\vec{y}.
\end{equation}
This system is sometimes called \emph{normal equation} in optimization community. Computing matrix on the left side of \eqref{eq:normal-equation} is prohibitively expensive as it is usually very large and dense. Even more prohibitively expensive is solving the system with this matrix as it is $\mathcal{O}(n^3)$ algorithm. One common approach is to drop either $\frac{1}{2\beta}M$ or $KM^{-1}K^T$ in the equation and to solve system with the remaining matrix, thus, using it as a preconditioner (see, for instance, \cite{RDW10}). However, this preconditioner would perform badly if none of the terms is dominant. We propose a completely different approach. Matrices $M$ and $K$ in \eqref{eq:normal-equation} are multilevel quasiseparable (see Examples \ref{ex:Laplace} and \ref{ex:tensor_product}). Thus, we can compute the Schur complement matrix in the left hand side of \eqref{eq:normal-equation} in quasiseparable arithmetic using $\mathcal{O}(n)$ arithmetic operations if quasiseparable orders stay small during computations. The last condition is true in practice because Schur complement is itself a discretized elliptic operator similar to Laplacian. In the previous example we have shown that it is more effective to use the proposed approach as a preconditioner rather than a direct solver. Thus, our approach to solving equation \eqref{eq:saddle-point2} consists of the following steps:
\begin{enumerate}
\item Invert matrix $M$ and form matrix $S=\frac{1}{2\beta}M+KM^{-1}K^T$ in quasiseparable matrices arithmetic.
\item Compute an approximate LDU decomposition of $S$ using Algorithm \ref{alg:LU} and some order reduction tolerance.
\item Use PCG method to solve $S\vec{\lambda}=\vec{y}$ with approximate LDU decomposition as a preconditioner.
\item Exploit computed $\vec{\lambda}$, $M^{-1}$ and the factorization \eqref{eq:LDU} to find $\vec{f}$ and $\vec{u}$.
\end{enumerate}

We have realized the proposed approach in practice. Tables \ref{tbl:saddle_point1} and \ref{tbl:saddle_point2} gather the numerical results we have obtained while solving the problem in Example \ref{ex:Poisson_control}.

\begin{table}
\centering
\caption{Time required to construct the normal equation matrix $S$ in the left hand side of \eqref{eq:normal-equation} with $\beta = 10^{-2}$, compute its approximate LDU decomposition and solve system with it to a specified tolerance using PCG method (number of iterations is given in brackets). Quasiseparable order used in order reduction algorithm was set to \textbf{1}.}
\label{tbl:saddle_point1}
\begin{tabular}[c]{rrrrr}
\multicolumn{5}{l}{$\beta = 10^{-2}$, $\mbox{tol}=10^{-4}$} \\
\midrule
$n$ & $S$ & LDU($S$) & PCG & total time \\
\midrule
12 & 0.73 & 1.12 & 0.0839 (3) & 1.93 \\
14 & 3.23 & 5.01 & 0.337 (3) & 8.57 \\
16 & 13.57 & 21.12 & 1.37 (3) & 36.05 \\
18 & 55.75 & 86.85 & 3.91 (2) & 146.51 \\
~ \\
\multicolumn{5}{l}{$\beta = 10^{-2}$, $\mbox{tol}=10^{-8}$} \\
\midrule
$n$ & $S$ & LDU($S$) & PCG & total time \\
\midrule
12 & 0.71 & 1.11 & 0.13 (5) & 1.95 \\
14 & 3.23 & 5.00 & 0.534 (5) & 8.76 \\
16 & 13.58 & 21.11 & 2.17 (5) & 36.86 \\
18 & 55.79 & 86.83 & 8.76 (5) & 151.38 \\
\bottomrule
\end{tabular}
\end{table}

\begin{table}
\centering
\caption{Time required to construct the normal equation matrix $S$ in the left hand side of \eqref{eq:normal-equation} with $\beta = 10^{-5}$, compute its approximate LDU decomposition and solve system with it to a specified tolerance using PCG method (number of iterations is given in brackets). Quasiseparable order used in order reduction algorithm was set to \textbf{1}.}
\label{tbl:saddle_point2}
\begin{tabular}[c]{rrrrr}
\multicolumn{5}{l}{$\beta = 10^{-5}$, $\mbox{tol}=10^{-4}$} \\
\midrule
$n$ & $S$ & LDU($S$) & PCG & total time \\
\midrule
12 & 0.71 & 1.12 & 0.0344 (1) & 1.86 \\
14 & 3.23 & 4.98 & 0.14 (1) & 8.35 \\
16 & 13.57 & 21.05 & 0.566 (1) & 35.19 \\
18 & 55.83 & 86.58 & 2.28 (1) & 144.70 \\
~ \\
\multicolumn{5}{l}{$\beta = 10^{-5}$, $\mbox{tol}=10^{-8}$} \\
\midrule
$n$ & $S$ & LDU($S$) & PCG & total time \\
\midrule
12 & 0.71 & 1.11 & 0.058 (2) & 1.88 \\
14 & 3.23 & 4.98 & 0.238 (2) & 8.45 \\
16 & 13.57 & 21.05 & 0.966 (2) & 35.59 \\
18 & 55.79 & 86.71 & 3.9 (2) & 146.40 \\
\bottomrule
\end{tabular}
\end{table}

\begin{figure}
\centering
\caption{Plots of the state, $u$, and the control, $f$, for $\beta=10^2$ and $10^{-2}$.}\label{fig:control}
\begin{tabular}{cc}
\includegraphics[scale = 0.35]{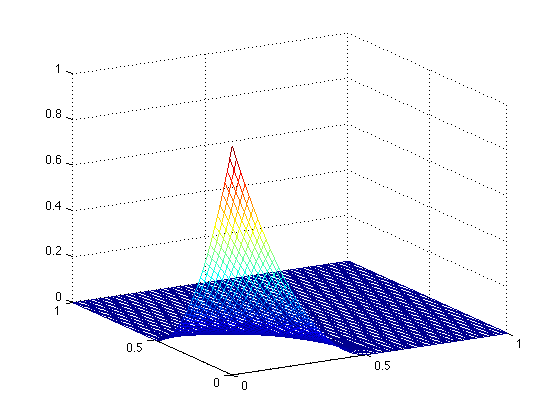}
&
\includegraphics[scale = 0.35]{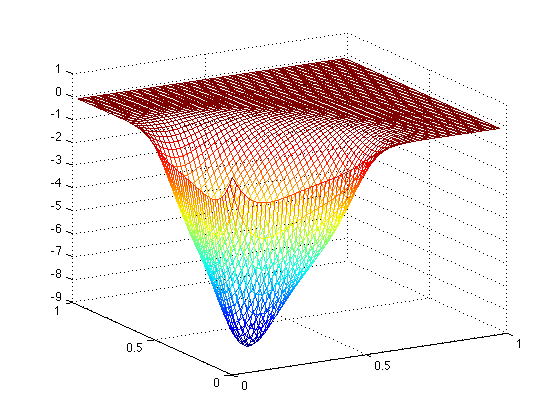}
\\
(a) $u$, $\beta = 10^2$, $\frac{\|u_h-\widehat{u}_h\|_2}{\|\widehat{u}_h\|_2}$=0.03
& 
(b) $f$, $\beta = 10^2$
\\ 
\includegraphics[scale = 0.35]{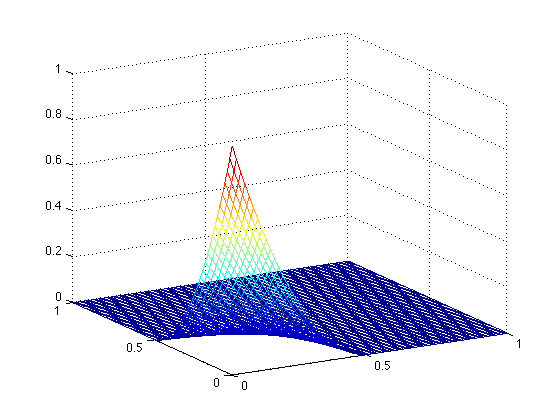}
&
\includegraphics[scale = 0.35]{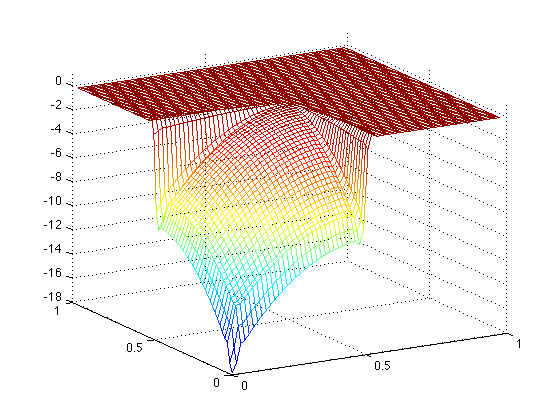}
\\
(c) $u$, $\beta = 10^{-2}$, $\frac{\|u_h-\widehat{u}_h\|_2}{\|\widehat{u}_h\|_2}$=1.9e-4
& 
(d) $f$, $\beta = 10^{-2}$
\end{tabular}
\end{figure}

Analyzing the results presented in Tables \ref{tbl:saddle_point1} and \ref{tbl:saddle_point2} we conclude that our preconditioner is mesh-independent and, as in the case of simple PDE problem, CPU time for the overall algorithm grows linearly with the size of the problem. Control function $f$ and obtained solution $u$ computed with the proposed algorithm for different values of the regularization parameter $\beta$ are presented in Figure \ref{fig:control}.

\section{Conclusion.}\label{sec:conclusion}
In this paper we have introduced a new class of rank-structured matrices called \emph{multilevel quasiseparable}. This matrices are low-parametric in the sense that only $\mathcal{O}(n)$ parameters are needed to store a matrix. Moreover, arithmetic operations and matrix decompositions can be performed in $\mathcal{O}(n)$ floating-point operations. Multilevel quasiseparable matrices extend the applicability of well-known quasiseparable matrices \cite{EG99a} from one-dimensional to multidimensional problems. In particular, we have shown that multilevel quasiseparable structure is well-suited to the description of discretized elliptic operators in 2D.

To demonstrate the usefulness of the new class of matrices we considered distributed control problem with a constraint in the form of a partial differential equation. Such problems were introduced by Lions and Mitter in \cite{L71}. In the course of solving them large-scale block matrices of saddle-point type arise. A straightforward way of solving these systems is to use block LU factorization. This is impractical as it requires direct inversion of large PDE matrices and further manipulations with them. However, we have shown that inverses of PDE matrices can be approximated by multilevel quasiseparable matrices with any desired accuracy and, hence, what was impossible in dense linear algebra became possible in structured linear algebra. 

Development of numerical methods for solving systems of saddle-point type is an important subject of modern numerical linear algebra, see an excellent survey paper by Benzi, Golub, and Liesen \cite{BGL05} for details. A large amount of work has been done on developing efficient preconditioners for such systems. In the current paper we have also proposed a new efficient mesh-independent preconditioner for saddle-point systems arising in PDE-constrained optimization. Performance of the new preconditioner is studied numerically.

There are several open questions left for future research. In particular, it is not clear what the range of applicability of multilevel quasiseparable matrices is. We have only considered problems with symmetric differential operators in our examples. However, theory does not require symmetry in the matrix and it may very well be that multilevel quasiseparable matrices are applicable to convection--diffusion and other unsymmetric operators as well. Next, robust techniques are needed to make the approach work on non-tensor grids. The biggest question though is how to extend order reduction Algorithm \ref{alg:reduction} to the block case to be able to use multilevel quasiseparable matrices in 3D and higher dimensions.


\bibliographystyle{plain}

\end{document}